\documentclass[a4paper,12pt]{amsart}

\numberwithin{equation}{section}
\theoremstyle{plain}
\newtheorem{theorem}{\sc \bf Theorem}[section]
\newtheorem{lemma}[theorem]{\sc \bf Lemma}
\newtheorem{corollary}[theorem]{\sc \bf Corollary}
\newtheorem{proposition}[theorem]{\sc \bf Proposition}
\newtheorem{claim}{\sc \bf Claim}
\theoremstyle{definition}
\newtheorem{definition}[theorem]{\sc \bf Definition}
\newtheorem{remark}[theorem]{\sc \bf Remark}
\newtheorem{example}[theorem]{\sc \bf Example}

\usepackage{amssymb}
\usepackage{amsfonts}
\usepackage{amsmath}
\usepackage{enumerate}

\usepackage{tikz}

\newcommand{\K}{{\mathbb{K}}}

\newcommand{\R}{{\mathbb{R}}}

\newcommand{\ali}[1]{\begin{align*}#1\end{align*}}
\newcommand{\alil}[1]{\begin{align}#1\end{align}}

\newcommand{\si}{\sigma}

\newcommand{\ph}{\varphi}
\newcommand{\phe}{\varphi(\epsilon,\cdot)}

\newcommand{\e}{\epsilon}
\newcommand{\p}{{\prime}}

\newcommand{\LR}{\mathcal{L}}
\let\MR=\relax
\newcommand{\MR}{\mathcal{M}}

\newcommand{\ite}[1]{\begin{enumerate}[(1)]#1\end{enumerate}}
\let\prop=\relax
\let\cor=\relax
\newcommand{\prop}[1]{\begin{proposition}#1\end{proposition}}
\newcommand{\props}{\begin{proposition}}
\newcommand{\prope}{\end{proposition}}
\newcommand{\cor}[1]{\begin{corollary}#1\end{corollary}}
\newcommand{\cors}{\begin{corollary}}
\newcommand{\core}{\end{corollary}}
\newcommand{\thm}[1]{\begin{theorem}#1\end{theorem}}
\newcommand{\thms}{\begin{theorem}}
\newcommand{\thme}{\end{theorem}}
\newcommand{\lem}[1]{\begin{lemma}#1\end{lemma}}
\newcommand{\lems}{\begin{lemma}}
\newcommand{\leme}{\end{lemma}}

\newcommand{\defis}{\begin{definition}}
\newcommand{\defie}{\end{definition}}

\newcommand{\exams}{\begin{example}}
\newcommand{\exame}{\end{example}}
\newcommand{\rem}[1]{\begin{remark}\normalfont #1\end{remark}}

\newcommand{\pros}{\begin{proof}}
\newcommand{\proe}{\end{proof}}
\newcommand{\prossq}{\begin{proof}}
\newcommand{\case}[1]{\begin{cases}#1\end{cases}}
\newcommand{\cd}{\cdot}
\newcommand{\up}{\upsilon}
\newcommand{\om}{\omega}
\newcommand{\ti}[1]{\tilde{#1}}

\newcommand{\diam}{{\mathrm{diam}}}

\newcommand{\qqqqquad}{\qquad\qquad\qquad}

\newcommand{\qqqqqqquad}{\qquad\qquad\qquad\qquad\qquad}

\newcommand{\MatIII}[1]{\left[\begin{array}{ccc}#1\end{array}\right]}

\newcounter{constants}
\makeatletter
\def\addconst{
\addtocounter{constants}{1}
\def\@currentlabel{\arabic{constants}}
\@currentlabel
}
\makeatother

\newcommand{\adl}[1]{\addconst\label{c:#1}}
\newcommand{\adr}[1]{\ref{c:#1}}

\newcommand{\dist}{{\mathrm{dist}}}

\begin{document}
\title[Dimension estimates in nonconformal GIFS]{Dimension estimates in nonconformal graph directed iterated function systems via asymptotic perturbation}
\author[H. Tanaka]{Haruyoshi Tanaka}
\address{
{\rm Haruyoshi Tanaka}\\
Department of Mathematics and Statistics\\
Wakayama Medical University\\
580, Mikazura, Wakayama-city, Wakayama, 641-0011, Japan
}
\email{htanaka@wakayama-med.ac.jp}
\begin{abstract}
We consider infinite graph-directed iterated function systems (GIFSs) whose contraction mappings are nonconformal. As our main result, we formulate asymptotic perturbations from conformal GIFSs to nonconformal GIFSs, and give the asymptotic behaviour of the Hausdorff dimension of the limit set of the perturbed system. We also investigate perturbed self-affine sets as special cases.
\end{abstract}
\maketitle
\section{Introduction}\label{sec:intro}
In this paper, we consider infinite graph-directed iterated function systems (GIFSs) whose contraction mappings are nonconformal. Let $E$ be a countable edge set of a directed graph. We introduce an $n$-order asymptotic perturbation $T_{e}(\e,x)$ of contraction mappings $T_{e}$ $(e\in E)$:
\ali
{
T_{e}(\e,\cd)=T_{e}+T_{e,1}\e+\cdots+T_{e,n}\e^{n}+o(\e^{n})
}
as $\e\to 0$ (see (\ref{eq:Tee=}) for detail). In particular, we assume that each $T_{e}$ is conformal but $T_{e}(\e,\cd)$ need not be conformal. Under suitable conditions, we give an $n$-order asymptotic expansion of the Hausdorff dimension of the limit set $J(\e)$ of the perturbed system $(T_{e}(\e,\cd))$:
\ali
{
\dim_{H}J(\e)=\dim_{H}J+s_{1}\e+\cdots+s_{n}\e^{n}+o(\e^{n})
}
as $\e \to 0$ (Theorem \ref{th:asymp_sol_quasiconformal}). This is a nonconformal version of the previous result \cite[Theorem 2.2]{T2023}. In particular, our potentials are more general than almost-additive potentials \cite{IY} (Remark \ref{rem:almostap}). Our result can be applied to affine maps case (Theorem \ref{th:asymp_affine}). Our result yields dimension estimates in affine maps not covered by previous results in \cite{CP,Falconer1992, Falconer, KR, Paulsen,Shmerkin,Solomyak} (see also Remark \ref{rem:affine1} and Remark \ref{rem:affine2}).
\smallskip
\par
To prove our main result, we use a technique of the solution of the equation of the pressure function $s\mapsto P(s\ph)$ for the physical potential $\ph$. If the functions $T_{e}$ are not conformal, then the Hausdorff dimension $\dim_{H}J$ of the limit set $J$ is not necessarily a solution of the pressure function $s\mapsto P(s\ph)$, but the number $P(((\dim_{H}J)/D)\ph)$ can be estimated by using the quasi-conformality constant $K$, where $D$ is the dimension of the Euclidean space (Corollary \ref{cor:GIFS_lower}). By using this fact and the previous result in \cite{T2023} for an asymptotic solution of the pressure function $s\mapsto P(s\ph(\epsilon,\cdot))$ for perturbed physical potential $\phe$, our main result is proved.
\smallskip
\par
In the next section \ref{sec:pre}, we introduce the notion of generalized quasiconformal maps and the definition of GIFS for state our main result. Our main theorem and this proof is shown in Section \ref{sec:main}. For a special case for nonconformal GIFS, we treat affine transformations and formulate an asymptotic perturbation of this in Section \ref{sec:app}. In particular, we demonstrate a concrete example of nonconformal mapping on $\R^{3}$ in Section \ref{sec:ex_aff}. To show our main result completely, we need upper and lower dimension estimates in nonconformal GIFS. These are stated in Section \ref{sec:aux}. In the appendix \ref{sec:asympsol}, we recall the result in \cite{T2023} for an asymptotic solution of the pressure function $s\mapsto P(s\ph(\epsilon,\cdot))$.
Moreover, we mention the mean valued inequalities on connected sets in Appendix \ref{sec:MVT}. In the final section \ref{sec:infnorm}, we will give fundamental properties of the map $\|\cd\|_{i}$ given by (\ref{eq:infnorm}).
\medskip
\\
\noindent
{\it Acknowledgment.}\ 
This study was supported by JSPS KAKENHI Grant Number 20K03636. \section{Preliminaries}\label{sec:pre}
\subsection{Generalized quasiconformal maps}
Fix integer $D\geq 1$. Let $O\subset \R^{D}$ be a domain. A mapping $f\,:\,O\to \R^{D}$ is a {\it generalized $K$-quasiregular} if $f$ is of $C^{1}$ and
there exists $1\leq K<\infty$ such that 
\alil
{
K^{-1}\|f^\p(x)\|^{D}\leq |\det f^\p(x)| \leq K (\|f^\p(x)\|_{i})^{D}\label{eq:def_quasireg}
}
for $x\in O$, where for a linear operator $\LR$ from a normed space $M$ to a normed space $N$, the number $\|\LR\|_{i}$ is defined by
\alil
{
\|\LR\|_{i}:=\inf_{x\in M\,:\,\|x\|=1}\|\LR x\|.\label{eq:infnorm}
}
Note that it is known that if the above condition replacing $|\det f^\p(x)|$ by $\det f^\p(x)$ holds, then $f$ is called $K$-quasiregular \cite{Rickman}.
A mapping $f\,:\, O\to f(O)$ is a {\it generalized quasiconfotmal} if  $f$ is a generalized $K$-quasiregular for some $K\geq 1$ and homeomorphism. A function $f\,:\,O\to \R^{D}$ is {\it conformal} if for any $x\in O$, there exists a constant $c_{\adl{con}}(x)>0$ such that $|f^\p(x)(u-v)|=c_{\adr{con}}(x)|u-v|$ for any $u,v\in \R^{D}$. Then we notice that $f$ is conformal if and only if $\|f^\p(x)\|=\|f^\p(x)\|_{i}$ for any $x\in O$ if and only if $f$ is a generalized $K$-quasiregular with $K=1$.
\subsection{Nonconformal graph directed function systems}\label{sec:AGIFS}
We introduce the following set $(G,(J_{v}),(O_{v}),(T_{e}))$ satisfying the conditions (G.1)-(G.4):
\begin{itemize}
\item[(G.1)] $G=(V,E,i(\cd),t(\cd))$ is a directed multigraph endowed with finite vertex set $V$, countable edge set $E$, and two maps $i(\cd)$ and $t(\cd)$ from $E$ to $V$. For each $e\in E$, $i(e)$ is called the initial vertex of $e$ and $t(e)$ called the terminal vertex of $e$.
\item[(G.2)] For each $v\in V$, $J_{v}$ is a compact and connected subset of $\R^{D}$ satisfying that the interior $\mathrm{int} J_{v}$ of $J_{v}$ is not empty, and $\mathrm{int} J_{v}$ and $\mathrm{int}J_{v^\p}$ are disjoint for $v^\p\in V$ with $v\neq v^\p$.
\item[(G.3)] For each $v\in V$, $O_{v}$ is a bounded, open and connected subset of $\R^{D}$ containing $J_{v}$.
\item[(G.4)] For each $e\in E$, $T_{e}$ is a $C^{1}$-diffeomorphism map from $O_{t(e)}$ to an open subset of $O_{i(e)}$ with $T_{e}(J_{t(e)})\subset J_{i(e)}$. Moreover, there exists $r\in (0,1)$ such that $\sup_{e\in E}\sup_{x\in O_{t(e)}}\|T_{e}^\p(x)\|\leq r$, where $\|T_{e}^\p(x)\|$ means the operator norm of $T_{e}^\p(x)$. 
\end{itemize}
For convenience, we call the set $(G,(J_{v}),(O_{v}),(T_{e}))$ satisfying (G.1)-(G.4) a graph iterated function system (GIFS for short). The {\it incidence matrix} $A$ of $G$ is a zero-one matrix indexed by $E$ such that $A(ee^\p)=1$ if $t(e)=i(e^\p)$ and $A(ee^\p)=0$ if $t(e)\neq i(e^\p)$. The code space is defined by
\alil
{
E^{\infty}=\{\om\in \prod_{n=0}^{\infty}E\,:\,A(\om_{n}\om_{n+1})=1 \text{ for all }n\geq 0\}.\label{eq:E^inf=}
}
Note that this set is also called the topological Markov shift with countable state space $E$ and with transition matrix $A$. The shift transformation $\si\,:\,E^{\infty}\to E^{\infty}$ is defined as $(\si\om)_{n}=\om_{n+1}$ for any $n\geq 0$. For $\theta\in (0,1)$, a metric $d_{\theta}$ on $E^{\infty}$ is given by $d_{\theta}(\om,\up)=\theta^{\min\{n\geq 0\,:\,\om_{n}\neq \up_{n}\}}$ if $\om\neq \up$ and $d_{\theta}(\om,\up)=0$ if $\om=\up$. Then $(E^{\infty},d_{\theta})$ is a complete separated metric space. The matrix $A$ is said to be {\it finitely irreducible} if there exists a finite subset $F$ of $\bigcup_{n=1}^{\infty}E^{n}$ such that for any $e,e^\p\in E$, $ewe^\p$ is a path on the graph $G$ for some $w\in F$. The coding map $\pi\,:\,E^{\infty}\to \R^{D}$ is well defined by
$\{\pi\om\}=\bigcap_{k=0}^{\infty}T_{\om_{0}\cdots\om_{k}}(J_{t(\om_{k})})$, where $T_{\om_{0}\cdots \om_{k}}$ means $T_{\om_{0}}\circ \cdots\circ T_{\om_{k}}$. The {\it limit set} of the system is given by the image $J:=\pi(E^{\infty})$.
A function $f\,:\,E^{\infty}\to \R$ is {\it acceptable} if there exists a constant $c_{\adl{acce}}\geq 1$ such that for any $e\in E$ and $\om,\up \in [e]$, $e^{f(\om)-f(\up)}\leq c_{\adr{acce}}$. For a function $f\,:\,E^{\infty}\to \R$, $P(f)$ means the {\it topological pressure} of $f$ which is defined by
\alil
{
P(f)=\lim_{n\to \infty}\frac{1}{n}\log \sum_{w\in E^{n}\,:\,[w]\neq \emptyset}\exp(\sup_{\om\in [w]}S_{n}f(\om)),\label{eq:P(f)=}
}
where $S_{n}f(\om):=\sum_{k=0}^{n-1}f(\si^{k}\om)$ and $[w]$ is the cylinder set $\{\om\in E^{\infty}\,:\,\om_{0}\cdots \om_{n-1}=w\}$ of the word $w\in E^{n}$. Note that if $f$ is acceptable then $P(f)$ exists in $[-\infty,+\infty]$ (see \cite{MU}).
We put
\alil
{
\overline{\ph}(\om)=&\log \|T_{\om_{0}}^{\p}(\pi\si\om)\|\label{eq:ophi}\\
\underline{\ph}(\om)=&\log(\|T_{\om_{0}}^{\p}(\pi\si\om)\|_{i})\label{eq:uphi}
}
for $\om\in E^{\infty}$.
We also impose the following conditions if necessary.
\ite
{
\item[$(G.5)_{C}$] (Conformally) Each $T_{e}$ is conformal.
\item[$(G.5)_{QC}$] There exists $K\geq 1$ such that each $T_{e}$ is generalized $K$-quasiregular.
\item[$(G.6)_{S}$] (Strong separated condition (SSC)) For $e,e^\p\in E$ with $e\neq e^\p$, $T_{e}(J_{t(e)}\cap J)\cap T_{e^\p}(J_{t(e^\p)}\cap J)=\emptyset$.
\item[$(G.6)_{O}$] (Open set condition (OSC)) For $e,e^\p\in E$ with $e\neq e^\p$, $T_{e}(\mathrm{int}J_{t(e)})\cap T_{e^\p}(\mathrm{int}J_{t(e^\p)})=\emptyset$.
\item[$(G.7)_{B}$] (Bounded distortion) There exist constants $c_{\adl{bd1}}>0$ and $0<\beta\leq 1$ such that for any $e\in E$ and $x,y\in O_{t(e)}$, $|\|T_{e}^\p(x)\|-\|T_{e}^\p(y)\||\leq c_{\adr{bd1}}\|T_{e}^\p(x)\| |x-y|^{\beta}$.
\item[$(G.7)_{SB}$] (Strongly bounded distortion) The condition $(G.7)_{B}$ is satisfied. Moreover, there exist constants $c_{\adl{bd2}}>0$ and $0<\beta\leq 1$ such that for any $e\in E$ and $x,y\in O_{t(e)}$, $|\|T_{e}^\p(x)^{-1}\|-\|T_{e}^\p(y)^{-1}\||\leq c_{\adr{bd2}}\|T_{e}^\p(x)^{-1}\| |x-y|^{\beta}$.
}
Note that $\|T^\p_{e}(x)^{-1}\|=(\|T^\p_{e}(x)\|_{i})^{-1}$ by Proposition \ref{prop:prop_infnorm}. If $(G.5)_{C}$ is satisfied, then $\|T_{e}^\p(x)\|=\|T_{e}^\p(x)\|_{i}$ and therefore the conditions $(G.7)_{B}$ and $(G.7)_{SB}$ are identical. Moreover, $\overline{\ph}=\underline{\ph}=:\ph$ is satisfied in this case. We say that a conformal GIFS $(G,(J_{v}),(O_{v}),(T_{e}))$ is {\it strongly regular} if $0<P(s\ph)<+\infty$ for some $s\geq 0$.
\rem
{\label{rem:almostap}
Assume that $(G,(J_{v}),(O_{v}),(T_{e}))$ is a nonconformal GIFS under the conditions $(G.5)_{QC}$, $(G.6)_{S}$ and $(G.7)_{SB}$. Put $\ph_{n}(\om):=\log\sup_{x\in J_{t(\om_{n-1})}}\|(T_{\om_{0}\cdots \om_{n-1}})^{\p}(x)\|$ for $\om\in E^{\infty}$ and $n\geq 1$, where $T_{\om_{0}\cdots \om_{n-1}}$ means $T_{\om_{0}}\circ\cdots\circ T_{\om_{n-1}}$. Then it is not hard to see that
\ali
{
-c_{\adr{kaap}}-(n+m)\log K+\ph_{n}(\om)+\ph_{m}(\si^{n}\om)\leq& \ph_{n+m}(\om)\\
\leq& c_{\adr{kaap}}+(n+m)\log K+\ph_{n}(\om)+\ph_{m}(\si^{n}\om)
}
for any $n,m\geq 1$ and $\om\in E^{\infty}$ for some constant $c_{\adl{kaap}}>0$.
Therefore, the sequence $\{\ph_{n}\}$ is not almost-additive potential \cite{IY} in general. Here the sequence of potentials $\{\psi_{n}\}$ is said to be {\it almost-additive potential} if there exist a constant $c_{\adl{aap}}>0$ such that for any $n,m\geq 1$ and $\om\in E^{\infty}$, $-c_{\adr{aap}}+\psi_{n}(\om)+\psi_{m}(\si^{n}\om)\leq \psi_{n+m}(\om)\leq c_{\adr{aap}}+\psi_{n}(\om)+\psi_{m}(\si^{n}\om)$.
}
\section{Main result and its proof}\label{sec:main}
We formulate an asymptotic perturbation of graph iterated function systems which is introduced in Section \ref{sec:AGIFS}.
Fix an integer $n\geq 0$ and a number $0<\beta\leq 1$. We assume the following:
\ite
{
\item[$(C.1)_{n}$] A set $(G,(J_{v}),(O_{v}),(T_{e}))$ is a strongly regular GIFS with the conformality condition $(G.5)_{C}$, the open set condition $(G.6)_{O}$ and the bounded distortion condition $(G.7)_{B}$. Assume also that the incidence matrix of the graph $G$ is finitely irreducible and the function $T_{e}$ is of class $C^{1+n+\beta}$ for each $e\in E$.
\item[$(C.2)_{n}$] Under the same graph $G$ and the sets $(J_{v})$ and $(O_{v})$ given in $(C.1)_{n}$, for each $\e\in(0,1)$, the set $(G,(J_{v}),(O_{v}),(T_{e}(\e,\cd)))$ is a GIFS with the quasi-conformality condition $(G.5)_{QC}$, the strongly separated condition $(G.6)_{S}$ and the strongly bounded distortion condition $(G.7)_{SB}$. Assume also the following (i)-(iii):
\ite
{
\item[(i)] For each $e\in E$, there exist functions $T_{e,k}\,:\,O_{t(e)}\to \R^{D}$ of class $C^{1+n-k+\beta}$ ($k=1,2,\dots, n$) and $\ti{T}_{e,n}(\e,\cd)\,:\,O_{t(e)}\to\R^{D}$ of class $C^{1+\beta(\e)}$ with $0<\beta(\e)\leq 1$ such that the function $T_{e}(\e,\cd)$ has the $n$-ordered asymptotic expansion
\alil
{
T_{e}(\e,\cd)=T_{e}+T_{e,1}\e+\cdots+T_{e,n}\e^{n}+\ti{T}_{e,n}(\e,\cd)\e^{n}\ \text{ on }J_{t(e)},\label{eq:Tee=}
}
and $\sup_{e\in E}\sup_{x\in J_{t(e)}}|\ti{T}_{e,n}(\e,x)|\to 0$ and $\sup_{e\in E}\sup_{x\in J_{t(e)}}\|\frac{\partial}{\partial x}\ti{T}_{e,n}(\e,x)\|\to 0$ as $\e\to 0$.
\item[(ii)] There exist constants $0<t(l,k)\leq 1$ ($l=0,1,\dots, n$,\ $k=1,\dots, n-l+1$) and $0<\ti{t}\leq 1$ such that (a) the function $x\mapsto T_{e,l}^{(k)}(x)/\|T_{e}^{\p}(x)\|^{t(l,k)}$ is bounded, $\beta$-H\"older continuous and its H\"older constant is bounded uniformly in $e\in E$, (b) $c_{\adl{Gbd3}}(\e):=\sup_{e\in E}\sup_{x\in J_{t(e)}}(\|\frac{\partial}{\partial x}\ti{T}_{e,n}(\e,x)\|/\|T_{e}^\p(x)\|^{\ti{t}_{0}})\to 0$ as $\e\to 0$, and (c) $\dim_{H}J/D>p(n)$. Here we put 
\alil
{
p(n):=
\case
{
\underline{p}/\ti{t},&n=0\\
\max\big\{\underline{p}+n(1-t_{1}), \underline{p}+\frac{n}{2}(1-t_{2}), \cdots, \underline{p}+\frac{n}{n}(1-t_{n}),\\
\qqqqqqquad\frac{\underline{p}}{t_{1}},\ \frac{\underline{p}}{t_{2}},\ \cdots, \frac{\underline{p}}{t_{n}},\ \underline{p}+1-\ti{t},\ \frac{\underline{p}}{\ti{t}}\big\},& n\geq 1.
}\label{eq:pn=}
}
\alil
{
t_{k}:=&\min\{\frac{1}{D}\sum_{p=1}^{D}t(i_{p},j_{p}+1)\,:\,i:=\sum_{q=1}^{D}i_{q}\text{ and }j:=\sum_{q=1}^{D}j_{q}\text{ satisfy}\label{eq:tk=}\\
&\qqqqquad i=k \text{ and }j=0 \text{ or }0\leq i<k \text{ and }1\leq j\leq k-i\}\nonumber\\
\ti{t}:=&\min\left\{t_{n},\ \ti{t}_{0},\ \frac{\ti{t}_{0}}{D}+\frac{D-1}{D}t(1,1),\dots,\ \frac{\ti{t}_{0}}{D}+\frac{D-1}{D}t(n,1)\right\}\label{eq:tt=}\\
\underline{p}:=&\inf\{s\geq 0\,:\,P(s\ph)<+\infty\}/D \text{ with }\ph(\om):=\log\|T_{\om_{0}}^\p(\pi\si\om)\|.\label{eq:p-=}
}
\item[(iii)] For any $e\in E$, the function $T_{e}(\e,\cd)$ is a generalized $K(\e)$-quasiregular and $K(\e)$ satisfies $K(\e)=1+\ti{K}_{n}(\e)\e^{n}$ with $\ti{K}_{n}(\e)\to 0$.
}
}
Here $P(s\ph)$ in (\ref{eq:p-=}) means the topological pressure of $s\ph$ defined by (\ref{eq:P(f)=}). Such a perturbed system was firstly given in \cite[Section 2.1]{T2023} under the case $K(\e)\equiv 1$. Now we state about the condition (ii). Roughly speaking, the condition (ii) includes the case where the Hausdorff dimension $\dim_{H}J(\e)$ of the limit set $J(\e)$ can be expanded to length $n$, but not to length $n+1$, namely $\dim_{H}J(\e)$ is not series at $\e=0$. On the other hand, the condition (ii) might be seemed a difficult to check. We will give in Section \ref{sec:affine} a reasonable assumption for the condition (ii).
Remark also that this condition (ii) automatically follows if there exists a finite subset $E_{0}\subset E$ such that $T_{e}(\e,\cd)$ does not depend on $\e$ for all $e\in E\setminus E_{0}$.
\smallskip
\par
Now we are in a position to state our main result:
\thms
\label{th:asymp_sol_quasiconformal}
Assume that the conditions $(C.1)_{n}$ and $(C.2)_{n}$ are satisfied. Then the Hausdorff dimension $\dim_{H}J(\e)$ of the limit set $J(\e)$ of the perturbed system $(G,(J_{v}),(O_{v}),(T_{e}(\e,\cd)))$ has the $n$-order asymptotic expansion $\dim_{H}J(\e)=\dim_{H}J+s_{1}\e+\cdots+s_{n}\e^{n}+o(\e^{n})$ as $\e \to 0$. In particular, each coefficient $s_{k}$ is given by using $T_{e}, T_{e,1},\dots, T_{e,k}$.
\thme
In order to show this theorem, we begin with the following lemma. 
Denoted by $\pi(\e,\cd)$ the coding map of $J(\e)$ for $\e>0$.
\lem
{[{\cite{T2016,T2023}}]\label{lem:asymp_pi}
Assume that the conditions $(C.1)_{n}$ and $(C.2)_{n}$ are satisfied. Choose any $r_{1}\in (r,1)$. Then there exist bounded $d_{r_{1}}$-Lipschitz continuous functions $\pi_{1},\pi_{2},\dots, \pi_{n}$ from $E^{\infty}$ to $\R^{D}$ and bounded continuous function $\ti{\pi}(\e,\cd)$ from $E^{\infty}$ to $\R^{D}$ such that $\pi(\e,\cd)=\pi+\pi_{1}\e+\cdots+\pi_{n}\e^{n}+\ti{\pi}_{n}(\e,\cd)\e^{n}$ and $\|\ti{\pi}_{n}(\e,\cd)\|_{\infty}:=\sup_{\om\in E^{\infty}}|\ti{\pi}_{n}(\e,\om)|\to 0$ as $\e\to 0$.
}
\pros
This assertion is guaranteed by the proof of \cite[Lemma 3.12]{T2023} without changes. Note that the proof of the lemma 3.12 in \cite{T2023} mostly depends on \cite[Lemma 3.1]{T2016}.
\proe
We put
\alil
{
\overline{\ph}(\e,\om)=&\textstyle\log \|\frac{\partial}{\partial x}T_{\om_{0}}(\e,\pi(\e,\si\om))\|\label{eq:ophie}\\
\underline{\ph}(\e,\om)=&\textstyle\log \|\frac{\partial}{\partial x}T_{\om_{0}}(\e,\pi(\e,\si\om))\|_{i}.\label{eq:uphie}
}
\lem
{\label{lem:conv_dim_nonconf}
Assume that the conditions $(C.1)_{n}$ and $(C.2)_{n}$ are satisfied. Then for any small $\e>0$, there exist a unique solution $s=\overline{s}(\e)\geq 0$ for $P(s\overline{\ph}(\e,\cd))=0$ and a unique solution $s=\underline{s}(\e)\geq 0$ for $P(s\underline{\ph}(\e,\cd))=0$ such that $\underline{s}(\e)\leq \dim_{H}J(\e)\leq \overline{s}(\e)$. Moreover, $\underline{s}(\e)$ and $\overline{s}(\e)$ converge to $\dim_{H}J$ both and hence $\dim_{H}J(\e)\to \dim_{H}J$.
}
\pros
First we show $\overline{s}(\e)\to \dim_{H}J/D$ as $\e\to 0$. We will check the conditions (g.1)-(g.5) in Appendix \ref{sec:asympsol} with $n=0$ for the functions $g(\e,\om):=\|\frac{\partial}{\partial x}T_{\om_{0}}(\e,\pi(\e,\si\om))\|$ and $g(\om):=\|T_{\om_{0}}^\p(\pi\si\om)\|$. The condition (g.1) with $n=0$ follows from the condition (g.5). By $\|T^\p_{e}(x)\|<1$, we see the condition (g.2).  The condition (g.3) and the condition $(G.7)_{B}$ for $(T_{e})$ are identical. When $n=0$, we ignore the condition (g.4). Therefore it is enough to prove the condition (g.5). By the condition $(C.2)_{n}$-(iii), we have
\ali
{
|\ti{g}_{0}(\e,\om)|=&|g(\e,\om)-g(\om)|\\
\leq& \|\textstyle\frac{\partial}{\partial x}T_{\om_{0}}(\e,\pi(\e,\si\om))-T_{\om_{0}}^\p(\pi\si\om)\|\\
\leq& \|\textstyle\frac{\partial}{\partial x}T_{\om_{0}}(\e,\pi(\e,\si\om))-T_{\om_{0}}^\p(\pi(\e,\si\om))\|+\|T_{\om_{0}}^\p(\pi(\e,\si\om))-T_{\om_{0}}^\p(\pi\si\om)\|\\
\leq& c_{\adr{Gbd3}}(\e)\|T_{\om_{0}}^\p(\pi(\e,\si\om))\|^{\ti{t}_{0}}+c_{\adr{Tbl}}\|T_{\om_{0}}^\p(\pi\si\om)\|^{t(0,1)}|\pi(\e,\si\om)-\pi\si\om|^{\beta}\\
\leq&c_{\adr{Gbd3}}(\e)((c_{\adr{bd1}}\|\pi(\e,\cd)-\pi\|^{\beta}+1)^{\ti{t}_{0}}\|T_{\om_{0}}^\p(\pi\si\om)\|^{\ti{t}_{0}}+c_{\adr{Tbl}}\|T_{\om_{0}}^\p(\pi\si\om)\|^{t(0,1)}\|\pi(\e,\cd)-\pi)\|^{\beta}\\
\leq&(c_{\adr{Gbd3}}(\e)(c_{\adr{bd1}}\|\ti{\pi}_{0}(\e,\cd)\|^{\beta}+1)^{\ti{t}_{0}}+c_{\adr{Tbl}}\|\ti{\pi}_{0}(\e,\cd)\|^{\beta})\|T_{\om_{0}}^\p(\pi\si\om)\|^{\ti{t}}=c_{\adr{ceg}}(\e)|g(\om)|^{\ti{t}},
}
where $c_{\adl{Tbl}}>0$ is a constant and $c_{\adl{ceg}}(\e)$ converges to $0$. Note that this convergence does not need the conformality of $T_{e}(\e,\cd)$.
Then (g.5) is guaranteed. $s(0):=\dim_{H}$ satisfies $s(0)>p(0):=\underline{p}/\ti{t}$ and $P(s(0)\log|g|)=0$ by a Bowen formula. By Theorem \ref{th:asympsol_Beq_Mgene}, there exists a unique solution $s=s(\e)$ for the equation $P(s\log|g(\e,\cd)|)=0$ such that $s(\e)\to s(0)$. The number $s(\e)$ equals $\overline{s}(\e)$ by the definition. Thus we obtain $\overline{s}(\e)\to \dim_{H}J$.
\smallskip
\par
On the other hand, we show $\underline{s}(\e)\to \dim_{H}J$. To do this, we will confirm the condition (g.5) for the functions $\underline{g}(\e,\om):=\|\frac{\partial}{\partial x}T_{\om_{0}}(\e,\pi(\e,\si\om))\|_{i}$ and the same $g(\om)=\|T_{\om_{0}}^\p(\pi\si\om)\|$. By using the condition $(C.2)_{n}$-$(iii)$, we notice
\ali
{
K(\e)^{-2/D}g(\e,\om)\leq \underline{g}(\e,\om)\leq g(\e,\om).
}
By $K(\e)^{-2/D}=1-(2/D)(1+\alpha \ti{K}_{n}(\e)\e^{n})^{-2/D}\ti{K}_{n}(\e)\e^{n}=1-c_{\adr{cK}}(\e)\e^{n}$ with $c_{\adl{cK}}(\e):=1-K(\e)^{-2/D}/\e^{n}$. Then $0\leq c_{\adr{cK}}(\e)\leq (2/D)\ti{K}_{n}(\e)\to 0$ is satisfied. We have
\ali
{
\underline{g}(\e,\om)-g(\om)\leq g(\e,\om)-g(\om)\leq c_{\adr{ceg}}(\e)|g(\om)|^{\ti{t}}
}
and
\ali
{
g(\om)-\underline{g}(\e,\om)\leq& g(\om)-g(\e,\om)+c_{\adr{cK}}(\e)\e^{n}g(\e,\om)\\
=&(1+c_{\adr{cK}}(\e)\e^{n})|g(\e,\om)-g(\om)|+c_{\adr{cK}}(\e)\e^{n}|g(\om)|\\
\leq &((1+c_{\adr{cK}}(\e)\e^{n})c_{\adr{ceg}}(\e)+c_{\adr{cK}}(\e)\e^{n})|g(\om)|^{\ti{t}}.
}
Thus we obtain
\ali
{
|\underline{g}(\e,\om)-g(\om)|\leq c_{\adr{ceg2}}(\e)|g(\om)|^{\ti{t}}
}
with $c_{\adl{ceg2}}(\e):=(1+c_{\adr{cK}}(\e)\e^{n})c_{\adr{ceg}}(\e)+c_{\adr{cK}}(\e)\e^{n}\to 0$. By virtue of Theorem \ref{th:asympsol_Beq_Mgene}, we get $P(\underline{s}(\e)\underline{g}(\e,\cd))=0$ and $\underline{s}(\e)$ converges to $\dim_{H}J$. 
\proe
\lems
[{\cite[Lemma 3.13]{T2023}}]\label{lem:asmyp_e^phe_MOREgene}
Assume that the conditions $(C.1)_{n}$ and $(C.2)_{n}$ are satisfied. Then the functions $g(\e,\om)=\det \frac{\partial}{\partial x}T_{\om_{0}}(\e,\pi(\e,\si\om))$ and $g(\om)=\det T_{\om_{0}}^{\prime}(\pi\si\om)$ satisfy the conditions (g.1)-(g.5) in Appendix \ref{sec:asympsol}.
\leme
\pros
Remark that the lemma 3.13 in \cite{T2023} is still satisfied without the conformality of $T_{\e}(\e,\cd)$.
\proe
\pros
[Proof of Theorem \ref{th:asymp_sol_quasiconformal}]
Let $g(\e,\om)=\det\frac{\partial}{\partial x}T_{\om_{0}}(\e,\pi(\e,\si\om))$ and  $g(\om)=\det T_{\om_{0}}^\p(\pi\si\om)$. Put $s(0)=\dim_{H}J/D$ and $\ph=\log|g|$. Then we have $P(s(0)\ph)=0$.  Observe that the condition $(C.2)_{n}$-$(iii)(c)$ yields $s(0)>p(n)$. Choose any compact neighborhood $I\subset (p(n),+\infty)$ of $s(0)$. By Theorem \ref{th:asympsol_Beq_Mgene}, for any small $\e>0$ there is a unique solution $s=s(\e)\in I$ for the equation $P(s\log|g(\e,\cd)|)=0$ and $s(\e)$ has the $n$-order asymptotic expansion (\ref{eq:t(e)=t0+t1e+...}):
\alil
{
s(\e)=s(0)+s_{1}\e+\cdots+s_{n}\e^{n}+\ti{s}_{n}(\e)\e^{n}.\label{eq:se=s0+s1e+...}
}
and $\ti{s}_{n}(\e)\to 0$. On the other hand, let $t(\e)=\dim_{H}J(\e)/D$. By using Corollary \ref{lem:conv_dim_nonconf}, $t(\e)\in I$ for any small $\e>0$. Moreover, Corollary \ref{cor:GIFS_lower} implies $|P(t(\e)\log|g(\e,\cd))|\leq t(\e)\log K(\e)$. Let $c_{\adl{cpt}}(\e)=P(t(\e)\log|g(\e,\cd))/\e^{n}$. We see
\ali
{
|c_{\adr{cpt}}(\e)|\leq |t(\e)\ti{K}_{n}(\e)|\leq (\max I)\ti{K}_{n}(\e)\to 0.
}
Now we notice the equation
\alil
{
P(t(\e)\log|g(\e,\cd)|-c_{\adr{cpt}}(\e)\e^{n})=0\label{eq:P(telogge-ce)=0}
}
since $-c_{\adr{cpt}}(\e)\e^{n}$ is constant.
We put 
\ali
{\psi\equiv 1,\ \psi_{1}=\cdots=\psi_{n}=0 \text{ and }\ti{\psi}_{n}(\e,\cd)\equiv (e^{-c_{\adr{cpt}}(\e)\e^{n}}-1)/\e^{n}.
}
The conditions $(\psi.1)$-$(\psi.3)$ in Appendix \ref{sec:asympsol} are immediately satisfied. To see the condition $(\psi.4)$, we have the estimate
\ali
{
|\ti{\psi}_{n}(\e,\om)|=e^{-\alpha c_{\adr{cpt}}(\e)\e^{n}}c_{\adr{cpt}}(\e)\leq e^{(\max I)\ti{K}_{n}(\e)\e^{n}}(\max I)\ti{K}_{n}(\e)=c_{\adr{ctpn}}(\e)|\psi(\om)|
}
with $c_{\adl{ctpn}}(\e):=e^{(\max I)\ti{K}_{n}(\e)\e^{n}}(\max I)\ti{K}_{n}(\e)\to 0$, where the first equation follows from Mean Valued Theorem with a number $\alpha\in [0,1]$. Note that (\ref{eq:P(telogge-ce)=0}) is equivalent to
\ali
{
P(t(\e)\log|g(\e,\cd)|+\log\psi(\e,\cd))=0
}
letting $\psi(\e,\cd):=\psi+\psi_{1}\e+\cdots+\psi_{n}\e^{n}+\ti{\psi}_{n}(\e,\cd)\e^{n}$.
By virtue of Theorem \ref{th:asympsol_Beq_Mgene} again, $t(\e)=\dim_{H}J(\e)/D$ has the $n$-order asymptotic expansion:
\ali
{
\dim_{H}J(\e)=Ds(0)+Ds_{1}\e+\cdots+Ds_{n}\e^{n}+\ti{t}_{n}(\e)\e^{n}
}
and $\ti{t}_{n}(\e)\to 0$. Here each coefficient $s_{k}$ is the same in (\ref{eq:se=s0+s1e+...}), namely each $s_{k}$ is decided by $T_{e}, T_{e,1},\dots, T_{e,k}$. Hence we obtain the assertion by noting $s(0)=\dim_{H}J/D$.
\proe
\section{Applications}\label{sec:app}
\subsection{Affine transformations}\label{sec:affine}
A mapping $T\,:\,\R^{D}\to \R^{D}$ is an {\it affine transformation} on $\R^{D}$ if $T$ has the form $T(x)=M(x)+a$ for some non-singular linear transformation $M$ acting on $\R^{D}$ and some vector $a\in \R^{D}$. 
\thm
{\label{th:GDMS_affine}
Let $(G,(J_{v}),(O_{v}),(T_{e}))$ be a GIFS with the finitely irreducible incidence matrix and with the strongly separated condition $(G.6)_{S}$. Assume also that each $T_{e}(x)=M_{e}x+a_{e}$ is an affine transformation. Then the Hausdorff dimension of the limit set $J$ satisfies $\underline{s}\leq \dim_{H}J\leq \overline{s}$, where $\underline{s}=\inf\{s\geq 0\,:\,P(s\underline{\ph})\leq 0\}$ and $\overline{s}=\inf\{s\geq 0\,:\,P(s\overline{\ph})\leq 0\}$ with $\underline{\ph}(\om)=\log\|M_{\om_{0}}\|$ and $\underline{\ph}(\om)=-\log\|M_{\om_{0}}^{-1}\|$.
}
\pros
It immediately follows that the condition $(G.7)_{SB}$ is guaranteed from $x\mapsto T^{\p}(x)=M_{e}$ is a constant function. Thus the assertion is valid by using Theorem \ref{th:GIFS_upper} and Theorem \ref{th:GIFS_lower}.
\proe
\rem
{\label{rem:affine1}
\item It is known in \cite{Falconer} that if the graph $G$ is a finite graph with singleton vertex and each $T_{e}x=M_{e}x+a_{e}$ is an affine transformation with $\|M_{e}\|<1/3$ and Lebesgue-almost all $a_{e}$, then the so-called Falconer dimension of $\{T_{e}\}$ is equal to the Hausdorff dimension of the self-affine set. In \cite{Solomyak}, this assertion is extended to the case $\|M_{e}\|<1/2$. On the other hand, \cite{Falconer1992, Paulsen} give the lower bounds for self-affine sets under the general case $\|M_{e}\|<1$ and under the strong separated condition (see also \cite{CP}). 
}
Now we formulate an asymptotic perturbation of affine maps. Fix an integer $n\geq 0$. We put $\ph(\om)=\log\|T_{\om_{0}}^\p(\pi\si\om)\|$, $\underline{s}=\{s\geq 0\,:\,P(s\ph)<+\infty\}$ and for $t\in (0,1]$
\ali
{
p_{n}(t):=\max(\underline{s}+Dn(1-t),\underline{s}/t).
}
Then we consider the following conditions:
\ite
{
\item[(F.1)] Let $(G,(J_{v}),(O_{v}),(T_{e}))$ be a GIFS with the finitely irreducible incidence matrix, with the open set condition $(G.7)_{O}$, and with strongly regular. Assume also that each $T_{e}(x)=M_{e}x+a_{e}$ is an affine transformation and $T_{e}$ is conformal, i.e. $M_{e}$ is similitude on $\R^{D}$.
\item[(F.2)]
\ite
{
\item[(i)] For each $\e>0$, a set $(G, (J_{v}), (O_{v}), (T_{e}(\e,\cd)))$ is a GIFS.
\item[(ii)] Each $T_{e}(\e,x)=M_{e}(\e)x+a_{e}(\e)$ is an affine transformation and there exist $t\in (0,1]$ with $p_{n}(t)<\dim_{H}J$, matrices $M_{e,k}$ and $\ti{M}_{e,n}(\e)$ and numbers $a_{e,k},\ti{a}_{e,n}(\e)\in \R^{D}$ $(1\leq k\leq n)$ such that
\ali
{
M_{e}(\e)=&M_{e}+M_{e,1}\e+\cdots+M_{e,n}\e^{n}+\ti{M}_{e,n}(\e)\e^{n}\\
a_{e}(\e)=&a_{e}+a_{e,1}\e+\cdots+a_{e,n}\e^{n}+\ti{a}_{e,n}(\e)\e^{n}
}
satisfying the finiteness $\max_{k}\sup_{e}\|M_{e,k}\|/\|M_{e}\|^{t}<\infty$, $\max_{k}\sup_{e}|a_{e,k}|<\infty$ and convergence $\sup_{e\in E}\|\ti{M}_{e,n}(\e)\|/\|M_{e}\|^{t}\to 0$, $\sup_{e\in E}|a_{e}(\e)|\to 0$ as $\e\to 0$.
\item[(iii)] Each $T_{e}(\e,x)$ is a generalized $K(\e)$-quasiregular with $K(\e)=1+o(\e^{n})$ as $\e\to 0$.
}
}
Then we obtain the following:
\thm
{\label{th:asymp_affine}
Assume that the conditions (F.1)(F.2) are satisfied. Then the Hausdorff dimension of the limit set of the set $(G, (J_{v}), (O_{v}), (T_{e}(\e,\cd)))$ has an $n$-order asymptotic expansion.
}
\pros
We put $T_{e,k}(x)=M_{e,k}x+a_{e,k}$ and $\ti{T}_{e,n}(\e,x)=\ti{M}_{e,n}(\e)x+\ti{a}_{e,n}(\e)$. We have the expansion $T_{e}(\e,\cd)=T_{e}+\sum_{k=1}^{n}T_{e,k}\e^{k}+\ti{T}_{e,n}(\e,\cd)\e^{n}$.
Since $T_{e,k}^\p=M_{e,k}$ and $T_{e,k}^{(i)}\equiv 0$ for all $i\geq 2$ hold, we can put $\ti{t}_{0}=t$, $t(k,1)=t$ for each $k=1,2,\dots, n$ and otherwise $t(k,l)\equiv 1$. 
Therefore $t_{k}=\ti{t}=t$ is satisfied. Thus we see $p(n)=\max(\underline{p}+n(1-t),\underline{p}/t)=\max(\underline{s}+n(1-t),\underline{s}/t)/D$ for any case $n\geq 0$. The assumption $p_{n}(t)>\dim_{H}J$ if and only if $p(n)>\dim_{H}J/D$. Hence the conditions $(C.1)_{n}$ and $(C.2)_{n}$ are fulfilled and the assertion is yielded from Theorem \ref{th:asymp_sol_quasiconformal}.
\proe
\rem
{\label{rem:affine2}
Note that \cite{KR} gave Falconer type formula for infinite affine IFS under contraction ratio $<1/2$. We stress that our above result does not need the restriction for contraction ratio. Remark also that \cite{Shmerkin} treated the continuity of the Hausdorff dimension of the self-affine set of finite affine IFS.
}
\subsection{A concrete example for affine transformation}\label{sec:ex_aff}
In this section, we give a concrete example for perturbed affine transformations. We define a map $T(\e,\cd)\,:\,\R^{3}\to\R^{3}$ by $T(\e,x)=M(\e)x+a(\e)$ with
\ali
{
M(\e)=
r\MatIII{
1/2&0&0\\
0&1/4&-\sqrt{3}/4\\
0&\sqrt{3}/4&1/4\\
}+
r\MatIII{
\e/4&0&0\\
0&\e/2&0\\
0&0&\e/2\\
}
}
and arbitrary choosing $a(\e)\in \R^{3}$. Note that $M(0)$ is the affine transformation that multiplies by $r$ in the $x$ direction, and rotates by $\pi/3$ degrees and multiplies by $r$ in the $yz$ direction. Namely, $M(\e)$ is a perturbation of $M(0)$. Let
\alil
{
K(\e)=2\sqrt{2}\frac{(\e^2+\e+1)^{3/4}}{(\e+2)^{3/2}}.\label{eq:Ke=}
}
\prop
{
The map $T(\e,\cd)$ is a generalized $K(\e)$-quasiregular map. In particular, $K(\e)=1+(9/16)\e^2+\cdots=1+o(\e)$ as $\e\to 0$.
}
\pros
We obtain
\ali
{
\|T(\e,x)\|=&\sup_{x\in \R^{3}\,:\,|x|=1}|M(\e)x|=\frac{1}{2}\sqrt{\e^2+\e+1}\\
\|T(\e,x)\|_{i}=&\inf_{x\in \R^{3}\,:\,|x|=1}|M(\e)x|=\frac{1}{4}(\e+2).
}
Hence $K(\e)=\sup_{x}(\|T(\e,x)\|/\|T(\e,x)\|_{i})^{3/2}$ equals (\ref{eq:Ke=}) and satisfies the assertion.
\proe
Consequently, if $(G,(J_{v}),(O_{v}),(T_{e}(\e,\cd)))$ is a GIFS such that $T_{e}(\e,\cd)\equiv T_{e}$ is conformal for $e\in E$ with $e\neq e_{0}$ and $T_{e_{0}}(\e,\cd)=T(\e,\cd)$, then the Hausdorff dimension of the limit set $J(\e)$ has a $1$-order asymptotic expansion $\dim_{H}J(\e)=\dim_{H}J+s_{1}\e+o(\e)$ as $\e\to 0$.
\section{Auxiliary results}\label{sec:aux}
In this section, we collect auxiliary results which are useful to show our main result. We start with the bounded distortion property. For convenience, we denote by $E^{*}$ the set of all finite path on the graph $G$.
\prop
{\label{prop:BD}
Assume that a set $(G,(J_{v}),(O_{v}),(T_{e}))$ is a GIFS and the condition $(G.7)_{B}$ is satisfied. Take a set of bounded open connected subsets $U_{v}$ of $O_{v}$ $(v\in V)$ so that $T_{e}U_{t(e)}\subset U_{i(e)}$ for any $e$. Then
\ite
{
\item there exists a constant $c_{\adl{BD}}\geq 1$ such that for any finite path $w=w_{1}\cdots w_{k}\in E^{*}$ and $x,y\in U_{t(w)}$,
\ali
{
\textstyle|\log(\prod_{i=1}^{k}\|T_{w_{i}}^\p(T_{w_{i+1}\cdots w_{k}}(x))\|)-\log(\prod_{i=1}^{k}\|T_{w_{i}}^\p(T_{w_{i+1}\cdots w_{k}}(y))\|)|\leq c_{\adr{BD}}|x-y|^{\beta};
}
\item if the condition $(G.7)_{SB}$ also holds, then exists a constant $c_{\adl{BD2}}\geq 1$ such that for any finite path $w=w_{1}\cdots w_{k}\in E^{*}$ and $x,y\in U_{t(w)}$,
\ali
{
\textstyle|\log(\prod_{i=1}^{k}\|T_{w_{i}}^\p(T_{w_{i+1}\cdots w_{k}}(x))^{-1}\|)-\log(\prod_{i=1}^{k}\|T_{w_{i}}^\p(T_{w_{i+1}\cdots w_{k}}(y))^{-1}\|)|\leq c_{\adr{BD2}}|x-y|^{\beta}.
}
}
}
\pros
(1) By virtue of the condition $(G.7)_{B}$, we notice $|T_{e}^\p(y)|\leq (1+c_{\adr{bd1}}|x-y|^{\beta})|T_{e}^\p(x)|$.
Put $x_{i}=T_{w_{i}\cdots w_{k}}(x)$ and $y_{i}=T_{w_{i}\cdots w_{k}}(y)$ and let $W_{v}=\bigcup_{z\in U_{v}}B(z,\delta)$ for small $\delta>0$ satisfying $W_{v}\subset O_{v}$. We also note the inequality $|x_{i}-y_{i}|\leq c_{\adr{MT}}r^{k-i+1}|x-y|$ from Proposition \ref{prop:T:mean} with the constant $c_{\adr{MT}}=c_{\adr{MT}}((U_{v}),(W_{v}))$. We obtain
\ali
{
&\Big|\prod_{i=1}^{k}\|T_{w_{i}}^\p(x_{i+1})\|-\prod_{i=1}^{k}\|T_{w_{i}}^\p(y_{i+1})\|\Big|\\
\leq&\sum_{j=1}^{k}\prod_{i=1}^{j-1}\|T_{w_{i}}^\p(x_{i+1})\|\left|\|T_{w_{j}}^\p(x_{j+1})\|-\|T_{w_{j}}^\p(y_{j+1})\|\right|\prod_{l=j+1}^{k}\|T_{w_{l}}^\p(y_{l+1})\|\\
\leq&\prod_{i=1}^{k}\|T_{w_{i}}^\p(x_{i+1})\|\sum_{j=1}^{k}c_{\adr{bd1}}|x_{j+1}-y_{j+1}|^\beta\prod_{l=j+1}^{k}(1+c_{\adr{bd1}}|x_{l+1}-y_{l+1}|^{\beta})\\
\leq&\prod_{i=1}^{k}\|T_{w_{i}}^\p(x_{i+1})\|\sum_{j=1}^{k}c_{\adr{bd1}}c_{\adr{MT}}^{\beta}r^{\beta(k-j)}|x-y|^\beta\prod_{l=j+1}^{k}(1+c_{\adr{bd1}}c_{\adr{MT}}^{\beta}r^{\beta(k-l)}|x-y|^{\beta})\\
\leq&c_{\adr{BD}}\prod_{i=1}^{k}\|T_{w_{i}}^\p(x_{i+1})\||x-y|^{\beta}
}
by putting $c_{\adr{BD}}=(c_{\adr{bd1}}c_{\adr{MT}}^{\beta}/(1-r^{\beta}))\prod_{i=0}^{\infty}(1+c_{\adr{bd1}}c_{\adr{MT}}^{\beta}\max_{v}(\diam U_{v})^{\beta})r^{i\beta})$, where the infinite product in $c_{\adr{BD}}$ is convergent by $\sum_{i=0}^{\infty}r^{i\beta}<+\infty$. Thus we get the assertion by using the basic inequality $|\log A-\log B|\leq |A-B|/\min(A,B)$ letting $A=\prod_{i=1}^{k}\|T_{w_{i}}^\p(x_{i+1})\|$ and $B=\prod_{i=1}^{k}\|T_{w_{i}}^\p(y_{i+1})\|$.
\smallskip
\\
(2) Together with Proposition \ref{prop:T:mean_inf}, a similar argument above implies the assertion
by putting $c_{\adr{BD2}}=(c_{\adr{bd2}}c_{\adr{iMT}}^{\beta}/(1-r^{\beta}))\prod_{i=0}^{\infty}(1+c_{\adr{bd2}}c_{\adr{iMT}}^{\beta}\max_{v}(\diam U_{v})^{\beta})r^{i\beta})$.
\proe
By virtue of the above proposition, we obtain the following distortion property.
\cor
{\label{cor:BD_v2}
Assume that a set $(G,(J_{v}),(O_{v}),(T_{e}))$ is a GIFS and the condition $(G.7)_{B}$ is satisfied. Take a set of bounded open connected subsets $U_{v}$ of $O_{v}$ $(v\in V)$ so that $T_{e}U_{t(e)}\subset U_{i(e)}$ for any $e$. Then
\ite
{
\item there exists a constant $c_{\adl{BDv}}\geq 1$ such that for any finite path $w=w_{1}\cdots w_{k}\in E^{*}$ and $x,y\in U_{t(w)}$, $c_{\adr{BDv}}^{-1}\leq \prod_{i=1}^{k}\|T_{w_{i}}^\p(T_{w_{i+1}\cdots w_{k}}(x))\|/\prod_{i=1}^{k}\|T_{w_{i}}^\p(T_{w_{i+1}\cdots w_{k}}(y))\|\leq c_{\adr{BDv}}$;
\item if $(G.7)_{SB}$ also holds, then there exists $c_{\adl{BDv2}}\geq 1$ such that for any finite path $w\in E^{*}$ and $x,y\in U_{t(w)}$, $c_{\adr{BDv2}}^{-1}\leq \prod_{i=1}^{k}\|T_{w_{i}}^\p(T_{w_{i+1}\cdots w_{k}}(x))\|_{i}/\prod_{i=1}^{k}\|T_{w_{i}}^\p(T_{w_{i+1}\cdots w_{k}}(y))\|_{i}\leq c_{\adr{BDv2}}$.
}
}
\pros
(1) By virtue of Proposition \ref{prop:BD}(1), we obtain the assertion by putting $c_{\adr{BDv}}=\exp(c_{\adr{BD}}\max_{v}(\diam U_{v})^{\beta})$. \\
(2) Letting $c_{\adr{BDv2}}=\exp(c_{\adr{BD2}}\max_{v}(\diam U_{v})^{\beta})$, we get the assertion.
\proe
\subsection{Upper dimension estimate in graph iterated function systems}\label{sec:Ues}
\thm
{\label{th:GIFS_upper}
Assume that a set $(G,(J_{v}),(O_{v}),(T_{e}))$ is a GIFS, the incidence matrix of $G$ is finitely irreducible and the condition $(G.7)_{B}$ is satisfied.
Then $\dim_{H}J\leq \overline{s}$, where $\overline{s}$ is given by $\overline{s}=\inf\{s\geq 0\,:\,P(s\overline{\ph})\leq 0\}$ and $\overline{\ph}$ is defined by (\ref{eq:ophi}).
}
\pros
Take $U_{v}=\bigcup_{z\in J_{v}}B(z,\delta)$ for small $\delta>0$ so that $U_{v}\subset O_{v}$ for each $v\in V$. Choose any $t>\inf\{s\geq 0\,:\,P(s\overline{\ph})\leq 0\}$. Then we see
\ali
{
0>P(t\overline{\ph})=&\lim_{n\to \infty}\frac{1}{n}\log\sum_{w\in E^{n}}\exp(\sup_{\om\in [w]}S_{n}(t\overline{\ph}))\\
=&\lim_{n\to \infty}\frac{1}{n}\log\sum_{w\in E^{n}}\sup_{\om\in [w]}\prod_{k=1}^{n}\|T_{w_{k}}^\p(\pi\si^{k+1}\om)\|^{t}.
}
Therefore there exists $n_{0}\geq 1$ such that for any $n\geq n_{0}$
\alil
{
\sum_{w\in E^{n}}\sup_{\om\in [w]}\prod_{k=1}^{n}\|T_{w_{k}}^\p(\pi\si^{k+1}\om)\|^{t}\leq \exp(nP(t\overline{\ph})/2).\label{eq:sumleq exp(nP(tph)/2)}
}
On the other hand, choose any $\delta>0$ and take $n\geq n_{0}$ so that $c_{\adr{MT}}r^{n}<\delta$. Since $\diam T_{w}J_{t(w)}\leq c_{\adr{MT}}\sup_{z\in U_{t(w)}}\|T_{w}^\p(z)\|\leq c_{\adr{MT}}r^{n}<\delta$ is satisfied, $(T_{w}J_{t(w)})_{w\in E^{n}}$ is a $\delta$-cover of $J$. In addition to bounded distortion property, we have
\ali
{
\mathcal{H}^{t}_{\delta}(J):=&\inf\{\sum_{k=1}^{\infty}(\diam C_{k})^{t}\,:\,\{C_{k}\} \text{ is a closed cover of } J \text{ and }\diam C_{k}< \delta\}\\
\leq& \sum_{w\in E^{n}}\diam(T_{w}J_{t(w)})^{t}\\
\leq &c_{\adr{MT}}^{t}\sum_{w\in E^{n}}\sup_{z\in O_{t(w)}}\|T_{w}^\p(z)\|^{t}\\
\leq&c_{\adr{MT}}^{t}\sum_{w\in E^{n}}\sup_{z\in U_{t(w_{i})}}\prod_{i=1}^{n}\|T_{w_{i}}^\p(T_{w_{i+1}\dots w_{n}}(z))\|^{t}\\
\leq&c_{\adr{MT}}^{t}c_{\adr{BDv}}^{t}\sum_{w\in E^{n}}\prod_{i=1}^{n}\|T_{w_{i}}^\p(T_{w_{i+1}\dots w_{n}}(\pi\si^{n+1}\om))\|^{t}\qquad (\because \text{Corollary }\ref{cor:BD_v2}(1))\\
\leq&c_{\adr{MT}}^{t}c_{\adr{BDv}}^{t}\exp(nP(t\overline{\ph})/2)\qquad (\because (\ref{eq:sumleq exp(nP(tph)/2)}))
}
for any $n\geq n_{0}$. Letting $n\to \infty$, $\exp(nP(t\overline{\ph})/2)\to 0$ and therefore we get $\mathcal{H}^{t}_{\delta}(J)=0$ for any $\delta>0$. Thus $\mathcal{H}^{t}(J)=\lim_{\delta\to 0}\mathcal{H}^{t}_{\delta}(J)=0$. This means $\dim_{H}J\leq t$. By arbitrarily choosing $t>\overline{s}=\inf\{s\geq 0\,:\,P(s\overline{\ph})\leq 0\}$, we obtain $\dim_{H}J\leq \overline{s}$. Hence the assertion is valid.
\proe
\subsection{Lower dimension estimate in graph iterated function systems}\label{sec:Les}
\thm
{\label{th:GIFS_lower}
Assume that a set $(G,(J_{v}),(O_{v}),(T_{e}))$ is a GIFS, the incidence matrix of $G$ is finitely irreducible and the conditions $(G.6)_{S}$ and $(G.7)_{SB}$ are satisfied. Then $\dim_{H}J\geq \underline{s}$, where $\underline{s}$ is given by $\underline{s}=\inf\{s\geq 0\,:\,P(s\underline{\ph})\leq 0\}$.
}
\prop
{\label{prop:prop_sepa}
Assume that a set $(G,(J_{v}),(O_{v}),(T_{e}))$ is a GIFS, the incidence matrix of $G$ is finitely irreducible and the conditions $(G.6)_{S}$ and $(G.7)_{SB}$ are satisfied. Take non-empty subsets $W_{v}\subset J_{v}$ $(v\in V)$ such that $T_{e}W_{t(e)}\subset W_{i(e)}$ for all $e\in E$. Put $\Delta(e,e^\p):=\dist(T_{e}W_{t(e)},T_{e^\p}W_{t(e^\p)})$ for $e,e^\p\in E$. Then there exists a constant $c_{\adl{sep}}>0$ such that for any $e,e^\p\in E$ with $e\neq e^\p$, and for any $w=w_{1}\cdots w_{k}\in E^{*}$ with $w\cd e,w\cd e^\p\in E^{*}$ and $z\in J_{t(w)}$, $\dist (W_{w\cd e},W_{w\cd e^\p})\geq c_{\adr{sep}}\Delta(e,e^\p)\prod_{i=1}^{k}\|T_{w_{i}}^\p(T_{w_{i+1}\cdots w_{k}}(z))\|_{i}$.
}
\pros
Take $U_{v}=B(J_{v},\delta)\subset O_{v}$ for small $\delta>0$ for all $v\in V$. Observe that $U_{v}$ is bounded open and connected, $T_{e}U_{t(e)}\subset U_{i(e)}$ and $W_{v}\subset U_{v}$. For $w\in E^{*}$, $e,e^\p\in E$ with $w\cd e,w\cd e^\p\in E^{*}$ and $e\neq e^\p$ and $z\in J_{t(w)}$, we have
\ali
{
\dist (W_{w\cd e},W_{w\cd e^\p})=&\inf_{x\in W_{t(e)},\ y\in W_{t(e^\p)}}|T_{w}(T_{e}(x))-T_{w}(T_{e^\p}(y))|\\
\geq&c_{\adr{iMT2}}\inf_{a\in U_{t(e)}}\|T_{w}^\p(a)\|_{i}\inf_{x\in W_{t(e)},\ y\in W_{t(e^\p)}}|T_{e}(x)-T_{e^\p}(y)|\\
\geq&c_{\adr{iMT2}}\inf_{a\in U_{t(e)}}\prod_{i=1}^{k}\|T_{w_{i}}^\p(T_{w_{i+1}\cdots w_{k}}(a))\|_{i}\Delta(e,e^\p)\\
\geq& c_{\adr{sep}}\prod_{i=1}^{k}\|T_{w_{i}}^\p(T_{w_{i+1}\cdots w_{k}}(z))\|_{i}\Delta(e,e^\p)
}
with $c_{\adl{iMT2}}=\inf_{e\in E}c_{\adr{iMT}}(W_{t(e)},U_{t(e)},W_{i(e)},U_{i(e)})$ and 
$c_{\adr{sep}}=c_{\adr{iMT2}}c_{\adr{BDv2}}^{-1}$ by using Corollary \ref{cor:BD_v2} and Proposition \ref{prop:T:mean_inf}. Hence the proof is complete.
\proe
Note that the following proposition imposes that the graph is finite.
\prop
{\label{prop:GIFS_lower_hole_ball}
Assume that a set $(G,(J_{v}),(O_{v}),(T_{e}))$ is a finite graph GIFS, $G$ is strongly connected and the conditions $(G.6)_{S}$ and $(G.7)_{SB}$ are satisfied. We put $W_{v}=\pi(\{\om\in E^{\infty}\,:\,i(\om)=v\})$ for $v\in V$. Notice that $J=\bigcup_{v\in V}W_{v}$, $T_{e}W_{t(e)}\subset W_{t(e)}$ and $W_{v}\subset J_{v}$. Then there exists a constant $c_{\adl{lowhb}}>0$, for any $\kappa\in (0,1)$, $k\geq 1$, $w\in E^{k}$, $z\in J_{t(w)}$ and $x\in T_{w}W_{t(w)}$, we have $J\cap B(x,c_{\adr{lowhb}}\kappa \prod_{i=1}^{k}\|T_{w_{i}}^\p(T_{w_{i+1}\cdots w_{k}}(z))\|_{i})\subset T_{w}W_{t(w)}$.
}
\pros
Choose any paths $w,\tau\in E^{k}$ with $w\neq \tau$. Take $1\leq l\leq k$ so that $\gamma:=w_{1}\cdots w_{l-1}=\tau_{1}\cdots \tau_{l-1}$ and $w_{l}\neq \tau_{l}$. By using Proposition \ref{prop:prop_sepa}, 
\ali
{
\dist(T_{w}W_{t(w)}, T_{\tau}W_{t(\tau)})\geq \dist(T_{\gamma\cd w_{l}}W_{t(w_{l})},T_{\gamma\cd \tau_{l}}W_{t(\tau_{l})})\geq&c_{\adr{lowhb}}\prod_{i=1}^{l-1}\|T_{w_{i}}^\p(T_{w_{i+1}\cdots w_{l-1}}(T_{w_{l}\cdots w_{k}}(z)))\|_{i}\\
\geq&c_{\adr{lowhb}}\prod_{i=1}^{k}\|T_{w_{i}}^\p(T_{w_{i+1}\cdots w_{k}}(z))\|_{i}\\
\geq &c_{\adr{lowhb}}\kappa\prod_{i=1}^{k}\|T_{w_{i}}^\p(T_{w_{i+1}\cdots w_{k}}(z))\|_{i}>0
}
by putting $c_{\adr{lowhb}}=c_{\adr{sep}}(\min_{e,e^\p\in E}\Delta(e,e^\p))$. By virtue of the condition $(G.6)_{S}$, we see $c_{\adr{lowhb}}>0$ from $\Delta(e,e^\p)>0$. Note that since $G$ is finite graph in this proposition, $E$ is a finite set. Thus we see
\ali
{
T_{\gamma}W_{t(\gamma)} \cap B(x,c_{\adr{lowhb}}\kappa\prod_{i=1}^{k}\|T_{w_{i}}^\p(T_{w_{i+1}\cdots w_{k}}(z))\|_{i})=\emptyset
}
for any $x\in T_{w}W_{t(w)}$. Hence the assertion holds by the equation $J=\bigcup_{w\in E^{k}}T_{w}W_{t(w)}$.
\proe
\pros[Proof of Theorem \ref{th:GIFS_lower}]
To prove this theorem, we consider the following two cases:\\
(Case 1): $G$ is finite graph. Recall the potential $\underline{\ph}$ defined in (\ref{eq:uphi}). Choose $\alpha>0$ so that $P(\alpha \underline{\ph})=0$. Denoted $\mu$ the Gibbs measure of $\alpha\underline{\ph}$ and put $\ti{\mu}=\pi\circ \pi^{-1}$. Let $A\subset K$ be a non-empty Borel subset, $x\in A$ and $\om=\pi^{-1}x$. Consider $\ti{\mu}(B(x,r))/r^{\alpha}$ for any small $r>0$. We let $\kappa:=(\min_{e\in E}\inf_{z\in J_{t(e)}}\|T_{e}^\p(z)\|_{i})^{-1}$. 
Choose any $0<r<\kappa^{-1}$. Set $r_{n}=\prod_{i=0}^{n}\|T_{\om_{i}}^\p(T_{\om_{i+1}\cdots \om_{n}}\pi\si^{n+1}\om)\|_{i}$ for each $k=0,1,\dots$. Then there exists a unique integer $k(r)>0$ such that $r_{k(r)}<r\leq r_{k(r)-1}$. 
since $r_{k(r)-1}\leq r_{k(r)}\kappa^{-1}$. Therefore we see $r_{k(r)}<r\leq \kappa^{-1}r_{k(r)}$. By virtue of Proposition \ref{prop:GIFS_lower_hole_ball},
\ali
{
&J\cap B(x,c_{\adr{lowhb}}\kappa \prod_{i=0}^{k(r)}\|T_{\om_{i}}^\p(T_{\om_{i+1}\cdots \om_{k(r)}}(\pi\si^{k(r)+1}\om))\|_{i})\subset T_{w}W_{t(w)}
\subset T_{\om_{0}\cdots \om_{k(r)}}(W_{t(\om_{k(r)})}),
}
where each set $W_{v}$ is given in Proposition \ref{prop:GIFS_lower_hole_ball}. Thus 
\ali
{
\ti{\mu}(B(x, c_{\adr{lowhb}}\kappa r))&\leq \ti{\nu}(B(x,c_{\adr{lowhb}}r_{k(r)}))\leq\ti{\mu}(T_{\om_{0}\cdots \om_{k(r)}}(K_{t(\om_{k(r)})}))
=\mu([\om_{0}\cdots \om_{k(r)}])\\
\leq &ce^{\alpha S_{k(r)+1}\underline{\ph}(\om)}
= c(\prod_{i=0}^{k(r)}\|T_{\om_{i}}^\p(\pi\si^{i+1}\om)\|_{i})^{\alpha}
= c (r_{n})^{\alpha}<c r^{\alpha}=c(c_{\adr{lowhb}}\kappa)^{-\alpha}(c_{\adr{lowhb}}\kappa r)^{\alpha},
}
where $c$ is a constant appearing in the definition of Gibbs measures.
By Frostman lemma, the $\alpha$-Hausdorff measure $\mathcal{H}^{\alpha}$ satisfies
\ali
{
\mathcal{H}^{\alpha}(K)\geq \frac{\ti{\nu}(K)}{c(c_{\adr{lowhb}}\kappa)^{-\alpha}}>0.
}
This means that $\alpha\leq \dim_{H}J$. We obtain the assertion under the finite graph.
\smallskip
\\
(Case II): $G$ is infinite graph.  Take a sequence of subgraphs $G_{n}=(V_{n},E_{n})$ of $G_{n}$ $(n\geq 1)$ with finite, strongly connected,  $V_{n}\subset V_{n+1}$, $E_{n}\subset E_{n+1}$, $\bigcup_{n}V_{n}=V$ and $\bigcup_{n}E_{n}=E$ (see \cite[Lemma 2.7.2]{MU} for the existence). Then $(G_{n},(J_{v})_{v\in V_{n}},(O_{v})_{v\in V_{n}}, (T_{e})_{e\in E_{n}})$ is a GIFS with finite graph and SSC.
We define $\underline{\ph}_{n}(\om):=\log\|T_{\om_{0}}^\p(\pi\si\om)\|_{i}$ for $\om\in E_{n}^{\infty}$. Take $\alpha_{n}>0$ satisfying $P(\alpha_{n}\underline{\ph}_{n})=0$.
Then $\dim_{H}J_{n}\geq \alpha_{n}$, where $J_{n}$ is the limit set of the subsystem $(G_{n},(J_{v})_{v\in V_{n}},(O_{v})_{v\in V_{n}}, (T_{e})_{e\in E_{n}})$. 
Since
\ali
{
0=P(\alpha_{n}\underline{\ph}_{n})\leq P(\alpha_{n}\underline{\ph}_{n+1}),
}
we get $\alpha_{n}\leq \alpha_{n+1}$. Let $\alpha_{\infty}:=\lim_{n\to \infty}\alpha_{n}$. Observe the inequality
\ali
{
\alpha_{\infty}\leq \lim_{n\to \infty}\dim_{H}J_{n}\leq \dim_{H}J.
}
For finite subset $E_{0}\subset E$, denoted by $G_{0}=(V,E_{0})$ the finite subgraph of $G$ and write $\underline{\ph}_{G_{0}}:=\underline{\ph}|_{(E_{0})^{\infty}}$. Since
\ali
{
P(\alpha_{\infty}\underline{\ph})=\sup_{G_{0}}P(\alpha_{\infty}\underline{\ph}_{G_{0}})\leq \sup_{n}P(\alpha_{\infty}\underline{\ph}_{n})
\leq \sup_{n}P(\alpha_{n}\underline{\ph}_{n})=0,
}
we see $\alpha_{\infty}\leq \underline{s}$, where the first equation is due to \cite[Theorem 2.1.5]{MU}.

On the other hand, we will show the converse $\underline{s}\geq \alpha_{\infty}$. Since $\alpha_{\infty}=0$ implies $\underline{s}=0\leq \dim_{H}J$, we may assume $\alpha_{\infty}>0$. Choose any $\delta>0$ so that $\alpha_{\infty}-\delta>0$. Then there exists $n_{0}\geq 1$ such that $\alpha_{\infty}-\delta<\alpha_{n_{0}}$ and
\ali
{
P((\alpha_{\infty}-\delta)\underline{\ph})\geq  P((\alpha_{\infty}-\delta)\underline{\ph}_{n_{0}})>P(\alpha_{n_{0}}\underline{\ph}_{n_{0}})=0.
}
We get $\underline{s}>\alpha_{\infty}-\delta$ and thus $\underline{s}\geq \alpha_{\infty}$. Hence $\underline{s}=\alpha_{\infty}\leq \dim_{H}J$.
\proe
The following is an easy corollary from Theorem \ref{th:GIFS_lower} and Theorem \ref{th:GIFS_upper}, but it is important role in our main theorem \ref{th:asymp_sol_quasiconformal}.
\cor
{\label{cor:GIFS_lower}
Under the same condition in Theorem \ref{th:GIFS_lower}, put $\ph(\om)=\log|\det(T_{\om_{0}}^\p(\pi\si\om))|$ and $p(0)=\inf\{s\geq 0\,:\,P(s\ph)<+\infty\}$. We also assume that the condition $(G.5)_{QC}$ is satisfied with the constant $K\geq 1$, and $p(0)<(\dim_{H}J)/D$. Then we have
\ali
{
|P(\frac{\dim_{H}J}{D}\ph)|\leq \frac{\dim_{H}J}{D}\log K.
}
}
\pros
Observe that the $K$-quasi-regularity (\ref{eq:def_quasireg}) implies that for any $s>0$
\ali
{
s\overline{\ph}-\frac{s}{D}\log K\leq \frac{s}{D}\ph(\om)\leq s\underline{\ph}+\frac{s}{D}\log K.
}
For any $s>\dim_{H}J$, we see $s>\underline{s}$ and $P(s\underline{\ph})<0$. Therefore
\ali
{
P(\frac{s}{D}\ph)\leq P(s\underline{\ph})+\frac{s}{D}\log K< \frac{s}{D}\log K.
}
Letting $s\to \dim_{H}J$, we obtain $P(((\dim_{H}J/D)\ph))\leq (\dim_{H}J/D)\log K$.

On the other hand, for any $s\in (s(0)D,\dim_{H}J)$, we have $s/D>s(0)$ and
\ali
{
P(\frac{s}{D}\ph)\geq& P(s\overline{\ph})-\frac{s}{D}\log K.
}
By $P((s/D)\ph)<+\infty$, $P(s\overline{\ph})$ is also finite. Moreover, the fact $s<\overline{s}$ yields $0<P(s\overline{\ph})<+\infty$. Thus
\ali
{
P(\frac{s}{D}\ph)>&-\frac{s}{D}\log K.
}
Letting $s\to \dim_{H}J$, we get $P((\dim_{H}J/D)\ph)\geq -(\dim_{H}J/D)\log K$. Hence the assertion holds.
\proe
\appendix
\section{Asymptotic solution of Bowen equation for perturbed potentials}\label{sec:asympsol}
We recall the result in \cite{T2023} which is the asymptotic solution of the equation of the pressure function $s\mapsto P(s\ph(\epsilon,\cdot))$ for perturbed potentials $\ph(\epsilon,\cdot)$ defined on the shift space with countable state space. To state our main result, we introduce some conditions for perturbed potentials. Let $G=(V,E,i(\cd),t(\cd))$ be a directed graph endowed with finite vertices $V$ and countable edges $E$. We use the notation $E^{\infty}$ and $d_{\theta}$ given in Section \ref{sec:AGIFS}. A function $f\,:\,E^{\infty}\to \K$ is called {\it weakly H\"older continuous} if the number $\sup_{e\in E}\sup_{\om,\up\in [e]\,:\,\om\neq \up}|f(\om)-f(\up)|/d_{\theta}(\om,\up)$ is finite for some $\theta\in (0,1)$. Note that `locally' might be used instead of `weakly'. Denoted by $\|f\|_{\infty}$ the supremum norm $\sup_{\om\in E^{\infty}}|f(\om)|$ of $f$. Let $n$ be a nonnegative integer. We consider the following conditions (g.1)-(g.5) for function $g(\e,\cd)\,:\,E^{\infty}\to \R$ with small parameter $\e\in (0,1)$:
\ite
{
\item[(g.1)] A function $g(\e,\cd)$ has the form
$g(\e,\cd)=g+g_{1}\e+\cdots+g_{n}\e^{n}+\ti{g}_{n}(\e,\cd)\e^{n}$
for some real-valued weakly H\"older continuous functions $g, g_{1},\dots, g_{n}, \ti{g}_{n}(\e,\cd)$ with $\|\ti{g}_{n}(\e,\cd)\|_{\infty}\to 0$ as $\e\to 0$.
\item[(g.2)] $g(\om)\neq 0$ for each $\om\in E^{\infty}$ and $\|g\|_{\infty}<1$.
\item[(g.3)] $|g(\om)-g(\up)|\leq c_{\adr{g3}}|g(\om)|d_{\theta}(\om,\up)$ for $\om,\up\in E^{\infty}$ with $\om_{0}=\up_{0}$ for some $c_{\adl{g3}}>0,\theta\in (0,1)$.
\item[(g.4)] $|g_{k}(\om)|\leq c_{\adr{g4}}|g(\om)|^{t_{k}}$ and $|g_{k}(\om)-g_{k}(\up)|\leq c_{\adr{g4_2}}|g(\om)|^{t_{k}}d_{\theta}(\om,\up)$ for any $\om,\up\in E^{\infty}$ with $\om_{0}=\up_{0}$ for some constants $c_{\adl{g4}},c_{\adl{g4_2}}>0$ and $0<t_{k}\leq 1$ for $k=1,2,\dots, n$.
\item[(g.5)] $|\ti{g}_{n}(\e,\om)|\leq c_{\adr{g5}}(\e)|g(\om)|^{\ti{t}}$ for any $\om\in E^{\infty}$ for some constants $0<\ti{t}\leq 1$ and $c_{\adl{g5}}(\e)>0$ with $c_{\adr{g5}}(\e)\to 0$.
}
Moreover, we assume that function $\psi(\e,\cd)\,:\,E^{\infty}\to \R$ satisfy the following conditions ($\psi.1$)-($\psi.4$):
\ite
{
\item[($\psi.1$)] A function $\psi(\e,\cd)$ has the form
$\psi(\e,\cd)=\psi+\psi_{1}\e+\cdots+\psi_{n}\e^{n}+\ti{\psi}_{n}(\e,\cd)\e^{n}$
for some real-valued weakly H\"older continuous functions $\psi, \psi_{1},\dots, \psi_{n}, \ti{\psi}_{n}(\e,\cd)$ with $\|\ti{\psi}_{n}(\e,\cd)\|_{\infty}\to 0$ as $\e\to 0$.
\item[($\psi.2$)] $\psi(\om)>0$ for any $\om\in E^{\infty}$.
\item[($\psi.3$)] $|\psi_{k}(\om)|\leq c_{\adr{psi3}}|\psi(\om)|$ and $|\psi_{k}(\om)-\psi_{k}(\up)|\leq c_{\adr{psi3_2}}|\psi(\om)|d_{\theta}(\om,\up)$ for any $\om,\up\in E^{\infty}$ with $\om_{0}=\up_{0}$ and for some $c_{\adl{psi3}},c_{\adl{psi3_2}}>0$ for $k=1,2,\dots, n$.
\item[($\psi.4$)] $|\ti{\psi}_{n}(\e,\om)|\leq c_{\adr{psi4}}(\e)|\psi(\om)|$ for any $\om\in E^{\infty}$ for some $c_{\adl{psi4}}(\e)>0$ with $c_{\adr{psi4}}(\e)\to 0$.
}
Let $\underline{p}=\inf\{p\geq 0\,:\,P(p\log |g|+\log\psi)<+\infty\}$ and recall the number $p(n)$ defined by (\ref{eq:pn=}). Then we obtained the following.
\thms
[{\cite[Theorem 1.1]{T2023}}]
\label{th:asympsol_Beq_Mgene}
Fix nonnegative integer $n$. Assume that the incidence matrix of $E^{\infty}$ is finitely irreducible and the conditions $(g.1)$-$(g.5)$ and $(\psi.1)$-$(\psi.4)$ are satisfied. Choose any $s(0)\in (p(n),+\infty)$ and any compact neighborhood $I\subset (p(n),+\infty)$ of $s(0)$. Let $p_{0}=P(s(0)\log|g|+\log\psi)$. Then there exist numbers $\e_{0}>0, s_{1},\dots, s_{n} \in \R$ such that the equation
$P(s\log|g(\e,\cd)|+\log\psi(\e,\cd))=p_{0}$
has a unique solution $s=s(\e)\in I$ for each $0<\e<\e_{0}$, and $s(\e)$ forms the asymptotic expansion
\alil
{
s(\e)=s(0)+s_{1}\e+\cdots+s_{n}\e^{n}+\ti{s}_{n}(\e)\e^{n}\label{eq:t(e)=t0+t1e+...}
}
and $\ti{s}_{n}(\e)\to 0$ as $\e\to 0$. In particular, each coefficient $s_{k}$ is defined by using $g,g_{1},\dots, g_{k}, \psi,\psi_{1},\dots,\psi_{k}$ (see \cite[Remark 3.9]{T2023}).
\thme
\section{Mean Valued inequality on connected sets}\label{sec:MVT}
\begin{proposition}[\cite{Patzschke}]\label{prop:T:mean}
Let $U$ be a non-empty open connected subset of $\R^{D}$ and $f\,:\,U\to Y$ a $C^{1}$ map from $U$ to a normed linear space $(Y,\|\cd\|_{Y})$. Take a bounded subset $V\subset U$ with $\dist(V,\partial U)>0$. Then
$\|f(x)-f(y)\|_{Y}\leq c_{\adr{MT}}\sup_{z\in U}\|f^{\prime}(z)\| |x-y|$ for each $x,y\in V$ with $c_{\adl{MT}}=c_{\adr{MT}}(V,U)=\max\{1,(\mathrm{diam}(V)/\mathrm{dist}(V,\partial U))\}$.
\end{proposition}
\prop
{[\cite{Patzschke}]\label{prop:T:mean_inf}
Let $O$ be a non-empty open subset of $\R^{D}$ and $f\,:\,O\to \R^{D}$ $C^{1}$-diffeomorphism. Take open connected subsets $U,U^\p\subset \R^{D}$ with $f(U)\subset U^\p$ and compact subsets $V\subset U, V^\p\subset U^\p$ with $f(V)\subset V^\p$. Then there exists $0<c_{\adl{iMT}}=c_{\adr{iMT}}(V,U,V^\p,U^\p)\leq 1$ such that $|f(x)-f(y)|\geq c_{\adr{iMT}}\inf_{z\in U}\|f^\p(z)^{-1}\|^{-1}|x-y|$ for any $x,y\in V$.
}
\section{Fundamental results of the map $\|\cd\|_{i}$}\label{sec:infnorm}
Let $M,N,L$ be normed spaces with $\dim M>0$. Recall that for a linear operator $\LR\,:\,M\to N$, the number $\|\LR\|_{i}$ is defined by $\inf_{x\in M\,:\,\|x\|=1}\|\LR x\|$ (the same as (\ref{eq:infnorm})).
\prop
{\label{prop:prop_infnorm}
Let $\LR\,:\,M\to N$ be a linear operator from a normed space $M$ with $\dim M\neq 0$ to a normed space $N$. Then
\ite
{
\item[(1)] $\|\LR\|_{i}=\inf_{x\in M\,:\,x\neq 0}\|\LR x\|/\|x\|=\inf_{x\in M\,:\,\|x\|\geq 1}\|\LR x\|$.
\item[(2)] If $\MR\,:\,N\to L$ is a linear operator, then $\|\MR\LR\|_{i}\geq \|\MR\|_{i}\|\LR\|_{i}$.
\item[(3)] If $\LR\,:\,M\to \LR(M)$ is injective, then $\|\LR\|_{i}=\|\LR^{-1}\|^{-1}$.
}
}
\pros
(1) We have
\ali
{
\inf_{x\in M\,:\,\|x\|\geq 1}\|\LR x\|\geq \inf_{x\in M\,:\,\|x\|\geq 1}\frac{\|\LR x\|}{\|x\|}\geq& \inf_{x\in M\,:\,x\neq 0}\left\|\LR \left(\frac{x}{\|x\|}\right)\right\|\\
\geq& \inf_{x\in M\,:\,\|x\|=1}\|\LR x\|=\|\LR\|_{i}\geq \inf_{x\in M\,:\,\|x\|\geq 1}\|\LR x\|.
}
(2) 
\ali
{
\|\MR\LR\|_{i}=\inf_{x\in M\,:\,x\neq 0}\|\MR(\LR x)\|&=\inf_{x\in M\,:\,x\neq 0,\ \LR x\neq 0}\|\MR(\LR x)\|\\
&\geq \inf_{y\in N\,:\,\|y\|=1}\|\MR y\|\inf_{x\in M\,:\,x\neq 0,\ \LR x\neq 0}\|\LR x\|\geq \|\MR\|_{i}\|\LR\|_{i}.
}
(3) We obtain
\ali
{
\|\LR^{-1}\|=\sup_{y\in \LR(M)\,:\,y\neq 0}\frac{\|\LR^{-1}y\|}{\|y\|}=\sup_{x\in M\,:\,x\neq 0}\frac{\|\LR^{-1}\LR x\|}{\|\LR x\|}=\left(\inf_{x\in M\,:\,x\neq 0}\frac{\|\LR x\|}{\|x\|}\right)^{-1}=(\|\LR\|_{i})^{-1}.
}
\proe

\endthebibliography
\end{document}